\documentstyle[twoside]{article}
\title{A note on an extension of Gelfond's constant}
\author{\sc Arjun K. Rathie$^a$ and Richard B. Paris$^b$\\
\\
${}^a\!$ Department of Mathematics, Vedant College of Engineering and Technology \\
(Rajasthan Technical University), Bundi, 323021, Rajasthan, India\\
E-Mail: arjunkumarrathie@gmail.com\\
${}^b\!$ Division of Computing and Mathematics,\\
Abertay University, Dundee DD1 1HG, UK\\
E-Mail: r.paris@abertay.ac.uk}
\begin{document}
\def\f#1#2{\mbox{${\textstyle \frac{#1}{#2}}$}}
\def\dfrac#1#2{\displaystyle{\frac{#1}{#2}}}
\def\boldal{\mbox{\boldmath $\alpha$}}
\newcommand{\bee}{\begin{equation}}
\newcommand{\ee}{\end{equation}}
\newcommand{\lam}{\lambda}
\newcommand{\ka}{\kappa}
\newcommand{\al}{\alpha}
\newcommand{\th}{\theta}
\newcommand{\om}{\omega}
\newcommand{\Om}{\Omega}
\newcommand{\fr}{\frac{1}{2}}
\newcommand{\fs}{\f{1}{2}}
\newcommand{\g}{\Gamma}
\newcommand{\br}{\biggr}
\newcommand{\bl}{\biggl}
\newcommand{\ra}{\rightarrow}
\renewcommand{\topfraction}{0.9}
\renewcommand{\bottomfraction}{0.9}
\renewcommand{\textfraction}{0.05}
\newcommand{\mcol}{\multicolumn}
\date{}
\maketitle
\pagestyle{myheadings}
\markboth{\hfill  {\it A.K. Rathie and R.B. Paris}  \hfill}
{\hfill {\it Gelfond's constant} \hfill}
\begin{abstract}
The aim of this note is to provide a natural extension of Gelfond's constant $e^\pi$ using a hypergeometric function approach. An extension is also found for the square root of this constant. A few interesting special cases are presented.
\vspace{0.4cm}

\setcounter{equation}{0}
\vspace{0.3cm}

\end{abstract}
\vspace{0.3cm}

\begin{center}
{\bf 1. \  Introduction}
\end{center}
\setcounter{section}{1}
\setcounter{equation}{0}
\renewcommand{\theequation}{\arabic{section}.\arabic{equation}}
In mathematics, Gelfond's constant, which is named after Aleksandr Gelfond, is given  by $e^\pi$. Like both $e$ and $\pi$, this constant is a transcendental number. The decimal expansion of Gelfond's constant is
\[e^\pi=23.1406926\,32779\ldots\ \]
and its continued fraction representation is given in \cite[A039661]{OES}. 

This number has a connection to the Ramanujan constant $e^{\pi\sqrt{163}}=(e^\pi)^{\sqrt{163}}$.
It is worth noting that this last number is almost an integer:
\[e^{\pi\sqrt{163}}\simeq 640320^3+744.\]
A geometrical occurrence of Gelfond's constant arises in the sum of even-dimension unit spheres with volume $V_{2n}=\pi^n/n!$. Then
\[\sum_{n=0}^\infty V_{2n}=e^\pi.\]

There are several ways of expressing Gelfond's constant, some of which are enumerated below:
\[e^\pi=(i^i)^{-2}\qquad (i=\sqrt{-1});\]
\[e^\pi=\left(\sum_{k=0}^\infty \frac{(-1)^k}{k!}\right)^{\!\!-4s},\qquad
s=\sum_{k=0}^\infty\frac{(-1)^k}{2k+1};\]
\[e^\pi=\bl(\prod_{k=1}^\infty k^{-\mu(k)/k}\br)^\sigma,\qquad \sigma=\sqrt{6\mbox{Li}_2(1)},\]
where $\mu(k)$ is the M\"obius function and $\mbox{Li}_n(x)$ is the polylogarithm function;
\[e^\pi={}_0F_1(-\!\!\!-;\fs;\pi^2/4)+\pi \,{}_0F_1(-\!\!\!-;\f{3}{2};\pi^2/4),\]
where ${}_0F_1(-\!\!\!-;a;z)$ is a generalised hypergeometric function that can be expressed in terms of  modified $I$-Bessel functions of order $\pm\fs$;
and finally
\bee\label{e11}
e^\pi={}_2F_1(i,-i; \fs;1)+2\,{}_2F_1(\fs+i,\fs-i; \f{3}{2};1),
\ee
where ${}_2F_1(a,b;c;z)$ is the well-known Gauss hypergeometric function \cite[p.~384]{DLMF}. The result (\ref{e11}) can be easily established by making use of the classical Gauss summation theorem
\bee\label{e12}{}_2F_1(a,b;c;1)=\frac{\g(c)\g(c-a-b)}{\g(c-a)\g(c-b)}
\ee
provided $\Re (c-a-b)>0$.

The natural extension of the summation theorem (\ref{e12}) to the ${}_3F_2$ hypergeometric series is available in the literature \cite{PBM}, which we shall write in the following manner:
\bee\label{e13}
{}_3F_2\bl(\begin{array}{c}a,b,d+1\\c+1,d\end{array}\!;1\br)=\frac{\g(c+1)\g(c-a-b)}{\g(c-a+1)\g(c-b+1)}\,\bl\{c-a-b+\frac{ab}{d}\br\}
\ee
provided $d\neq 0, -1, -2,\ldots$ and $\Re (c-a-b)>0$.
The aim of this note is to provide a natural extension of Gelfond's constant (\ref{e11}), and also its square root, with the help of the result (\ref{e13}). A few interesting results closely related to Gelfond's constant and its square root are also given.

\vspace{0.6cm}

\begin{center}
{\bf 2. \  Extension of Gelfond's constant}
\end{center}
\setcounter{section}{2}
\setcounter{equation}{0}
\renewcommand{\theequation}{\arabic{section}.\arabic{equation}}
The natural extension of Gelfond's constant to be established here is given in the following theorem.
\newtheorem{theorem}{Theorem}
\begin{theorem}
For $d_1, d_2\neq 0, -1, -2, \ldots\,$, the following result holds true:
\[e^\pi\bl(\frac{1}{5d_1}+\frac{15}{32d_2}+\frac{23}{80}\br)+e^{-\pi}\bl(\frac{1}{5d_1}-\frac{15}{32d_2}-\frac{7}{80}\br)\]
\bee\label{e21}
={}_3F_2\bl(\begin{array}{c}i,-i,d_1+1\\\f{3}{2},\,d_1\end{array}\!;1\br)+2\,{}_3F_2\bl(\begin{array}{c}\fs+i,\fs-i,d_2+1\\\ \f{5}{2},\,d_2\end{array}\!;1\br).
\ee
\end{theorem}

\noindent {\it Proof.}\ \ \ The derivation of (\ref{e21}) follows from application of the summation formula (\ref{e13}). We have
\[{}_3F_2\bl(\begin{array}{c}i,-i,d_1+1\\\f{3}{2},\,d_1\end{array}\!;1\br)=(e^\pi+e^{-\pi})\bl(\frac{1}{10}+\frac{1}{5d_1}\br)\]
and
\[{}_3F_2\bl(\begin{array}{c}\fs+i,\fs-i,d_2+1\\\ \f{5}{2},\,d_2\end{array}\!;1\br)=(e^\pi-e^{-\pi})\bl(\frac{3}{32}+\frac{15}{64d_2}\br).\]
Insertion of these summations into the right-hand side of (\ref{e21}) then yields the result asserted by the theorem $\Box$
\vspace{0.6cm} 

\begin{center}
{\bf 3. \  Corollaries}
\end{center}
\setcounter{section}{3}
\setcounter{equation}{0}
\renewcommand{\theequation}{\arabic{section}.\arabic{equation}}
In this section, we mention some interesting special cases of our main result in (\ref{e21}).
\medskip

\newtheorem{corollary}{Corollary}
\begin{corollary} \ \ In (\ref{e21}), if we take $d_1=2/(5n-1)$ and $d_2=15/(2(8n-3))$ for positive integer $n$, then we obtain after a little calculation the following result:
\bee\label{e31}
n e^\pi={}_3F_2\bl(\begin{array}{c}i,-i,\frac{5n+1}{5n-1}\\\f{3}{2},\,\frac{2}{5n-1}\end{array}\!;1\br)+2\,{}_3F_2\bl(\begin{array}{c}\fs+i,\fs-i,\frac{16n+9}{2(8n-3)}\\\f{5}{2},\,\frac{15}{2(8n-3)}\end{array}\!;1\br).
\ee
\end{corollary}
In particular, when $n=1$ we recover Gelfond's constant (\ref{e11}). For $n=2, 3$ we find respectively the following results related to (\ref{e11}):
\bee
2e^\pi={}_3F_2\bl(\begin{array}{c}i,-i,\frac{11}{9}\\\f{3}{2},\,\frac{2}{9}\end{array}\!;1\br)+2\,{}_3F_2\bl(\begin{array}{c}\fs+i,\fs-i,\,\frac{41}{26}\\\f{5}{2},\frac{15}{26}\end{array}\!;1\br)
\ee
and
\bee
3e^\pi={}_3F_2\bl(\begin{array}{c}i,-i,\frac{8}{7}\\\f{3}{2},\,\frac{1}{7}\end{array}\!;1\br)+2\,{}_3F_2\bl(\begin{array}{c}\fs+i,\fs-i,\frac{19}{14}\\\f{5}{2},\,\frac{5}{14}\end{array}\!;1\br).
\ee
\medskip

\begin{corollary} \ \ In (\ref{e21}), if we take $d_1=2/(5n-1)$ and $d_2=-15/(2(8n+3))$ for positive integer $n$, then we obtain after a little calculation the following result:
\bee\label{e32}
n e^{-\pi}={}_3F_2\bl(\begin{array}{c}i,-i,\frac{5n+1}{5n-1}\\\f{3}{2},\,\frac{2}{5n-1}\end{array}\!;1\br)+2\,{}_3F_2\bl(\begin{array}{c}\fs+i,\fs-i,\frac{16n-9}{2(8n+3)}\\\f{3}{2},\,\frac{-15}{2(8n+3)}\end{array}\!;1\br).
\ee
\end{corollary}
In particular, for $n=1, 2, 3$ we find respectively the following results:
\bee
e^{-\pi}= {}_2F_1\bl(\begin{array}{c}i,\,-i\\\f{1}{2}\end{array}\!;1\br)+2\,{}_3F_2\bl(\begin{array}{c}\fs+i\,\fs-i,\, \frac{7}{22}\\\f{3}{2},-\frac{15}{22}\end{array}\!;1\br),
\ee
\bee
2e^{-\pi}={}_3F_2\bl(\begin{array}{c}i,\,-i,\ \frac{11}{9}\\\f{3}{2},\,\frac{2}{9}\end{array}\!;1\br)+2\,{}_3F_2\bl(\begin{array}{c}\fs+i,\fs-i,\frac{23}{38}\\\f{3}{2},-\frac{15}{38}\end{array}\!;1\br)
\ee
and
\bee
3e^{-\pi}={}_3F_2\bl(\begin{array}{c}i,-i,\frac{8}{7}\\\f{3}{2},\,\frac{1}{7}\end{array}\!;1\br)+2\,{}_3F_2\bl(\begin{array}{c}\fs+i,\fs-i,\frac{13}{18}\\ \f{3}{2},-\frac{5}{18}\end{array}\!;1\br).
\ee
\medskip

\begin{corollary} \ \ In (\ref{e21}), if we take $d_1=1/(2(10n-1))$ and $d_2=-5/2$ for positive integer $n$, then we obtain after a little calculation the following result:
\bee\label{e33}
n(e^\pi+e^{-\pi}) ={}_3F_2\bl(\begin{array}{c}i,-i,\frac{20n-1}{2(10n-1)}\\\f{3}{2},\,\frac{1}{2(10n-1)}\end{array}\!;1\br)+2\,{}_3F_2\bl(\begin{array}{c}\fs+i,\fs-i,-\frac{3}{2}\\\f{3}{2},\,-\frac{5}{2}\end{array}\!;1\br).
\ee
\end{corollary}
In particular, for $n=1, 2, 3$ we find respectively the following results:
\bee
e^{\pi}+e^{-\pi}= {}_3F_2\bl(\begin{array}{c}i,-i, \frac{19}{18}\\\f{3}{2},\f{1}{18}\end{array}\!;1\br)+2\,{}_3F_2\bl(\begin{array}{c}\fs+i,\fs-i,-\frac{3}{2}\\\f{3}{2},\,-\frac{5}{2}\end{array}\!;1\br),
\ee
\bee
2(e^{\pi}+e^{-\pi})={}_3F_2\bl(\begin{array}{c}i,-i,\frac{39}{38}\\\f{3}{2},\,\frac{1}{38}\end{array}\!;1\br)+2\,{}_3F_2\bl(\begin{array}{c}\fs+i,\fs-i,-\frac{3}{2}\\\f{3}{2},\,-\frac{5}{2}\end{array}\!;1\br)
\ee
and
\bee
3(e^{\pi}+e^{-\pi})={}_3F_2\bl(\begin{array}{c}i,-i,\frac{59}{58}\\\f{3}{2},\,\frac{1}{58}\end{array}\!;1\br)+2\,{}_3F_2\bl(\begin{array}{c}\fs+i,\fs-i,-\frac{3}{2}\\\f{3}{2},\,-\frac{5}{2}\end{array}\!;1\br).
\ee
Similarly other results can be obtained.
\vspace{0.6cm}

\begin{center}
{\bf 4. \  The square root of Gelfond's constant: $e^{\pi/2}$}
\end{center}
\setcounter{section}{4}
\setcounter{equation}{0}
\renewcommand{\theequation}{\arabic{section}.\arabic{equation}}
Expressions for the square root of Gelfond's constant are:
\[e^{\pi/2}=i^{-i};\]
\[e^{\pi/2}=\left(\sum_{k=0}^\infty \frac{(-1)^k}{k!}\right)^{\!\!-2s},\qquad
s=\sum_{k=0}^\infty\frac{(-1)^k}{2k+1};\]
\bee\label{e41a}
e^{\pi/2}={}_2F_1\bl(\!\!\begin{array}{c} i, -i\\ \fs\end{array}\!;\frac{1}{2}\br)+\sqrt{2}\, {}_2F_1\bl(\!\!\begin{array}{c} \fs+i, \fs-i\\ \f{3}{2}\end{array}\!;\frac{1}{2}\br)
\ee
together with the inverse expression
\bee\label{e41b}
e^{-\pi/2}={}_2F_1\bl(\!\!\begin{array}{c} i, -i\\ \fs\end{array}\!;\frac{1}{2}\br)-\sqrt{2}\, {}_2F_1\bl(\!\!\begin{array}{c} \fs+i, \fs-i\\ \f{3}{2}\end{array}\!;\frac{1}{2}\br).
\ee
The results in (\ref{e41a}) and (\ref{e41b}) can be obtained by evaluating the first hypergeometric function by the second Gauss theorem and the second hypergeometric function by Bailey's theorem viz.
\[{}_2F_1\bl(\!\!\begin{array}{c} a, b\\ \fs (a+b+1)\end{array}\!;\frac{1}{2}\br)=\frac{\g(\fs) \g(\fs a+\fs b+\fs)}{\g(\fs a+\fs) \g(\fs b+\fs)},\]
\[{}_2F_1\bl(\!\!\begin{array}{c} a, 1-a\\ c\end{array}\!;\frac{1}{2}\br)=\frac{\g(\fs c) \g(\fs c+\fs)} {\g(\fs c+\fs a) \g(\fs c-\fs a+\fs)}.\]

We now derive the analogue of Theorem 1 by making use of the extension of the second Gauss and Bailey's theorems
applied to ${}_3F_2$ series. These are given by \cite{KRR}:
\[
{}_3F_2\bl(\!\!\begin{array}{c}a, b, d+1\\ \fs(a+b+3), d\end{array}\!;\frac{1}{2}\br)=\frac{\g(\fs)\g(\fs a+\fs b+\f{3}{2})\g(\fs a-\fs b-\fs)}{\g(\fs a-\fs b+\f{3}{2})}\]
\bee\label{e42}
\times\bl\{\frac{\fs(a+b-1)-ab/d}{\g(\fs a+\fs) \g(\fs b+\fs)}+\frac{(a+b+1)/d-2}{\g(\fs a) \g(\fs b)}\br\}
\ee
and
\[
{}_3F_2\bl(\!\!\begin{array}{c}a, 1-a, d+1\\ c+1, d\end{array}\!;\frac{1}{2}\br)=2^{-c} \g(\fs) \g(c+1)\]
\bee\label{e43}
\times\bl\{\frac{2/d}{\g(\fs c+\fs a) \g(\fs c-\fs a+\fs)}+\frac{1-(c/d)}{\g(\fs c+\fs a+\fs) \g(\fs c-\fs a+1)}\br\},
\ee
provided $d\neq 0, -1, -2, \ldots\ $.
Then we have the following theorem:
\begin{theorem}
For $d_1, d_2\neq 0, -1, -2, \ldots\,$, the following result holds true:
\[e^{\pi/2}\bl(\frac{1}{10d_1}+\frac{3}{16d_2}+\frac{27}{40}\br)+e^{-\pi/2}\bl(\frac{3}{10d_1}-\frac{21}{16d_2}+\frac{11}{40}\br)\]
\bee\label{e44}
={}_3F_2\bl(\begin{array}{c}i,-i,d_1+1\\\f{3}{2},\,d_1\end{array}\!;\frac{1}{2}\br)+\sqrt{2}\,{}_3F_2\bl(\begin{array}{c}\fs+i,\fs-i,d_2+1\\\ \f{5}{2},\,d_2\end{array}\!;\frac{1}{2}\br).\ee
\end{theorem}

\noindent {\it Proof.}\ \ \ In the first ${}_3F_2$ series use (\ref{e42}) and in the second ${}_3F_2$ series use (\ref{e43}) together with standard properties of the gamma function. $\Box$

\medskip

\begin{corollary}
If in (\ref{e44}) we take $d_1=1/(7n-5)$ and $d_2=15/(24n-14)$ for positive integer $n$ then we find
\bee\label{e45}
ne^{\pi/2}={}_3F_2\bl(\begin{array}{c}i,-i,\frac{7n-4}{7n-5}\\\f{3}{2},\,\frac{1}{7n-5}\end{array}\!;\frac{1}{2}\br)+\sqrt{2}\,{}_3F_2\bl(\begin{array}{c}\fs+i,\fs-i,\frac{24n+1}{24n-14}\\\ \f{5}{2},\,\frac{15}{24n-14}\end{array}\!;\frac{1}{2}\br).
\ee
\end{corollary}
When $n=1$ we recover (\ref{e41a}). For $n=2, 3$ we find respectively the following results:
\bee
2e^{\pi/2}={}_3F_2\bl(\begin{array}{c}i,-i,\frac{10}{9}\\\f{3}{2},\,\frac{1}{9}\end{array}\!;\frac{1}{2}\br)+\sqrt{2}\,{}_3F_2\bl(\begin{array}{c}\fs+i,\fs-i,\frac{49}{34}\\\ \f{5}{2},\,\frac{15}{34}\end{array}\!;\frac{1}{2}\br)
\ee
and
\bee
3e^{\pi/2}={}_3F_2\bl(\begin{array}{c}i,-i,\frac{17}{16}\\\f{3}{2},\,\frac{1}{16}\end{array}\!;\frac{1}{2}\br)+\sqrt{2}\,{}_3F_2\bl(\begin{array}{c}\fs+i,\fs-i,\frac{73}{58}\\\ \f{5}{2},\,\frac{15}{58}\end{array}\!;\frac{1}{2}\br).
\ee
Similarly other results can be obtained.

To conclude we note that an obvious extension of (\ref{e11}) is
\bee\label{e46}
e^{\pi\lambda}={}_2F_1(i\lambda,-i\lambda; \fs;1)+2\lambda \,{}_2F_1(\fs+i\lambda,\fs-i\lambda; \f{3}{2};1),
\ee
where $\lambda$ is real. This yields the alternative expression for $e^{\pm\pi/2}$ given by
\bee
e^{\pm\pi/2}={}_2F_1(\fs i,-\fs i; \fs; 1)\pm \,{}_2F_1(\fs+\fs i, \fs-\fs i; \f{3}{2}; 1).
\ee
In addition, if we choose $\lambda=\sqrt{19}$, $\sqrt{43}$, $\sqrt{67}$ and $\sqrt{167}$ in (\ref{e46}), we obtain hypergeometric function representations of the four largest Heegner numbers, viz.
\begin{eqnarray*}
e^{\pi\sqrt{19}}&=&8.85479\ 77768\ 01543\ 19497\ldots\times 10^5\label{e23}\\
&\doteq& 96^3+744-0.22;\nonumber\end{eqnarray*}
\begin{eqnarray*}
e^{\pi\sqrt{43}}&=&8.84736\ 74399\ 97774\ 66034\ldots \times 10^8\label{e24}\\
&\doteq& 960^3+744-0.00022;\nonumber\end{eqnarray*}
\begin{eqnarray*}
e^{\pi\sqrt{67}}&=&1.47197\ 95274\ 39999\,98662\ldots \times 10^{11}\label{e25}\\
&\doteq& 5280^3+744-0.0000013;\nonumber\end{eqnarray*}
\begin{eqnarray*}
e^{\pi\sqrt{163}}&=&2.62537\ 41264\ 07687\ 43999\ldots \times 10^{17}\label{e26}\\
&\doteq& 640320^3+744-7.5\times 10^{-13},\nonumber\end{eqnarray*}
the last number being Ramanujan's constant.

\end{document}